\documentclass [12pt,a4paper]{amsart}
\usepackage[utf8]{inputenc}
\usepackage[T1]{fontenc}
\usepackage{lmodern}
\usepackage {amssymb}
\usepackage {amsmath}
\usepackage {amsthm}
\usepackage{mathtools}
\usepackage{xfrac}
\usepackage{mathrsfs}
\usepackage[margin=1in]{geometry}

\usepackage{url}

\usepackage{enumerate}
\usepackage{hyperref}
\hypersetup{
    colorlinks=true,
    linkcolor=blue,
    filecolor=magenta,      
    urlcolor=cyan,
    citecolor=magenta
}
\usepackage{cleveref}

\DeclareMathOperator {\C} {\mathbb{C}}
\DeclareMathOperator {\R} {\mathbb{R}}

\DeclareMathOperator {\Z} {\mathbb{Z}}
\DeclareMathOperator {\N} {\mathbb{N}}
\DeclareMathOperator {\Q} {\mathbb{Q}}

\DeclareMathOperator {\dom} {\mathrm{dom}}

\newtheorem{thm}{Theorem}[section]
\newtheorem{cor}[thm]{Corollary}
\newtheorem{lem}[thm]{Lemma}
\newtheorem{prop}[thm]{Proposition}
\newtheorem{claim}[thm]{Claim}
\theoremstyle{definition}
\newtheorem*{defn}{Definition}

\newtheorem{ex}[thm]{Example}

\theoremstyle{remark}
\newtheorem{remark}[thm]{Remark}
\numberwithin{equation}{thm}

\title{Bialgebraic varieties of the Gamma function}

\author{Sebastian Eterovi\'c}
\email{sebastian.eterovic@univie.ac.at}
\address{Kurt G\"odel Research Center, Universit\"at Wien, 1090 Wien, Austria}

\author{Adele Padgett}
\email{adele.lee.padgett@univie.ac.at}
\address{Kurt G\"odel Research Center, Universit\"at Wien, 1090 Wien, Austria}

\author{Roy Zhao}
\email{rhzhao@caltech.edu}
\address{Department of Mathematics, California Institute of Technology, Pasadena, CA, USA}

\date{}

\thanks{SE was supported  by the Austrian Science Fund (FWF) 10.55776/ESP1584024.}



\begin{document}

\begin{abstract}
   We characterize the bialgebraic varieties of the $\Gamma$ function, that is, if $V,W\subseteq\C^n$ are irreducible affine algebraic varieties which satisfy $\dim V =\dim W$ and $\Gamma(V)\subseteq W$, then the equations defining $V$ (and hence also $W$) either give an equality between coordinates, or set some coordinates to be constant.
   We also classify the case where $V$ and $W$ have the same dimension as varieties over the field of algebraic functions.
\end{abstract}

\maketitle

\section{Introduction}

In this paper we study some of the algebraic functional properties of Euler's $\Gamma$ function.
We recall that the $\Gamma$ function is meromorphic on $\C$, having a simple pole at every non-positive integer, and satisfies $\Gamma(z+1) = z\Gamma(z)$.
We recall precise definitions in \S\ref{sec:prelims}.

The motivation behind this work is to make progress towards understanding the functional properties of $\Gamma$ in the spirit of functional transcendence results known in the literature as `Ax--Lindemann--Weierstrass' theorems (ALW). 
An ALW type theorem is a statement about a transcendental holomorphic function $f:U\to\C^n$, where $U\subseteq \C^n$ is some open complex domain, which describes the maximal algebraic varieties contained in $f^{-1}(W)$, where $W\subseteq\C^n$ is any algebraic variety.
Equivalently, such a theorem describes the algebraic varieties $V\subseteq\C^n$ such that $\dim\overline{f(V\cap U)}^{\mathrm{Zar}}<n$.
Our main results take a first step in this direction for the case where $f$ is taken to be the $n$-fold function $\Gamma:U_\Gamma^n\to\C^n$ defined as $(z_1,\ldots,z_n)\mapsto(\Gamma(z_1),\ldots,\Gamma(z_n))$, where $U_\Gamma\coloneqq \C\setminus\Z_{\leq0}$.
Since $\Gamma$ is a non-vanishing meromorphic function, one can equivalently work with the entire function $f=\frac{1}{\Gamma}$ if one does not want to worry about the poles.
We will often take this approach during the proof of our main theorem.

Functional transcendence results of this type (and even stronger results known as Ax--Schanuel theorems) have been proved for a wide array of functions: for the complex exponential function \cite{ax}, the modular $j$-function \cite{pila:andre-oort, pila-tsimerman:ax-schanuel}, more general modular functions \cite{cfn:alw}, and uniformization maps of Shimura varieties \cite{mpt}, among others \cite{bakker-tsimerman:asvhs,bcfn:ax-schanuel,cfn:alw,gao-klingler:ax-schanuel}. 

Throughout the paper, subvarieties of $\C^n$ will be defined by polynomials in the ring $\C[X_1,\ldots,X_n]$.
Given an algebraic variety $W$ and a holomorphic function $f:\C^n\to W$, we say that a subvariety $V$ of $\C^n$ is \emph{bialgebraic for $f$} if $\dim \overline{f(V)}^{\mathrm{Zar}} = \dim V$. 
Note that there is a family of varieties which are bialgebraic for the $n$-fold Cartesian product of any (set-theoretic) function.
\begin{defn}
    We say that an algebraic variety $V\subseteq\C^n$ is \emph{trivially bialgebraic} if $V$ is defined by equations of the following forms:
    \begin{enumerate}[(i)]
        \item $X_i=X_k$, for some $i,k\in\{1,\ldots,n\}$ distinct, 
        \item $X_i = w$, for some $i\in\{1,\ldots,n\}$ and some $w\in\C$.
    \end{enumerate}
\end{defn}

Our first main result shows that $\Gamma$ has no interesting bialgebraic varieties. 

\begin{thm}
\label{thm:equidimALW}
    Suppose $V,W \subset \C^n$ are irreducible algebraic varieties satisfying $\Gamma(V\cap U_\Gamma^n)\subseteq W$.
    Then $V$ is trivially bialgebraic if and only if $\dim V = \dim W$,.
\end{thm}

Our second main result shows that one can find interesting algebraic properties of $\Gamma$ if one works with the graph of $\Gamma$ instead of just the image of $\Gamma$. 

\begin{thm}
\label{thm:bialggraph}
    Suppose  $V\subseteq\C^n$ and $W\subseteq\C^{2n}$ are irreducible algebraic varieties satisfying that for all $(z_1,\ldots,z_n)\in V\cap U_\Gamma^n$ we have $(z_1,\ldots  ,z_n,\Gamma(z_1),\ldots,\Gamma(z_n)) \in W$. 
    Then $\dim W = 2\dim V$ if and only if there is a partition $\{1,\ldots,n\}=A\sqcup B\sqcup C$ such that
    \begin{enumerate}[(a)]
        \item the coordinate projection $\pi_A:V\to\C^A$ mapping to the coordinate with index in $A$ is dominant,
        \item the coordinate projection $\pi_A:V\to\C^C$ mapping to the coordinate with index in $A$ is a singleton and for every $i\in\{1,\ldots,n\}\setminus C$ the projection $\pi_i:V\to \C$ onto the $i$-th coordinate is dominant,
        \item the set $B$ has size at least $2$ and for all $i\in B$ there are $j\in B$ and $b_{ij}\in \Z$ such that $i\neq j$ and the relation $X_i =X_j+b_{ij}$ holds on $V$.
    \end{enumerate}
\end{thm}

Theorem \ref{thm:bialggraph} can also be interpreted as bialgebraicity result in the following sense: if $V\subseteq\C^n$ is a variety of dimension $d$ and $\alpha_1,\ldots,\alpha_n:U\to\C$ are algebraic functions defined on some open domain in $U\subseteq\C^d$ parameterizing $V$, then the theorem classifies those $V$ for which
\[\operatorname{tr.deg.}_{\C(\alpha_1,\ldots,\alpha_n)}\C(\alpha_1,\ldots,\alpha_n,\Gamma(\alpha_1),\ldots,\Gamma(\alpha_n))=\dim V.\]

\subsection{Functional properties and bialgebraic varieties}
\label{subsec:functrans}
For many interesting transcendental functions $f$ (such as those listed above: complex exponentiation, the modular $j$-function, etc.), the algebraic functional properties of $f$ are explained by its bialgebraic varieties.
That is to say, given an algebraic variety $W\subseteq\C^n$, the maximal algebraic varieties contained in $f^{-1}(W)$ are bialgebraic for $f$ when $f$ is one of those functions.  
However, this is not the case for $\Gamma$. 
One simple example is to let $V$ be the line in $\C^3$ parametrized by
\[V \coloneqq \{(z, z + 1, z + 2) : z \in \C\} \subset \C^3.\]
Observe that $V$ is not contained in a proper trivially bialgebraic variety.
Using the functional identities $\Gamma(z+2)=(z+1)\Gamma(z+1)=\left(\frac{\Gamma(z+1)}{\Gamma(z)}+1\right)\Gamma(z+1)$ one sees that for every $(a,b,c)\in V$ with $a\notin\Z_{\leq 0}$, the point $(\Gamma(a),\Gamma(b),\Gamma(c))$ belongs to the hypersurface $W\subseteq\C^3$ defined by the equation
\[X_2^2 + X_1X_2 = X_1X_3.\]
One may check (manually or using Theorem \ref{thm:equidimALW}) that in this case $\dim\overline{\Gamma(V\cap U_\Gamma^3)}^{\mathrm{Zar}} = 2$ and so $\overline{\Gamma(V\cap U_\Gamma^3)}^{\mathrm{Zar}} = W$. 

On the other hand, the dimension of the Zariski closure of the image of $V$ under the map 
\[
    (z_1,z_2,z_3) \mapsto (z_1,z_2,z_3,\Gamma(z_1),\Gamma(z_2),\Gamma(z_3))
\]
also has dimension 2, so this example satisfies the hypotheses of Theorem \ref{thm:bialggraph}.

\subsection{Related work in difference algebra}
The recent work of L.~Di Vizio and F.~Pellarin \cite{divizio-pellarin} also refers to a differential Ax--Lindemann--Weierstrass type statement for $\Gamma$, although it goes in a different direction than our motivation described earlier.
In \cite[Theorem 1]{divizio-pellarin}, the authors prove a strong result regarding the differential algebraic dependencies that occur among a set of functions of the form $\Gamma(z-f(z))$, where $f(z)$ is a function which is algebraic over the field of 1-periodic meromorphic functions on $\C$. 
Yet, their result does not consider the case where $f$ is algebraic over $\C(z)$, unless $f$ is constant.
Instead, our two main theorems and \cite[Theorem 1]{divizio-pellarin} complement each other by tackling different instances of functional transcendence questions for $\Gamma$, and are motivated by \cite[\S3]{EP25} and \cite[Question 2]{divizio-pellarin}.

\subsection{About the methods}
There are well-established methods for proving functional transcendence results for a certain class of transcendental functions.
For instance, if the map $f$ is obtained as the quotient map by the action of some lattice contained in an algebraic group acting algebraically on the domain of $f$, then the restriction of the function to a fundamental domain may be part of an o-minimal expansion of the reals, in which case one can use o-minimal complex analysis (possibly aided by point-counting) to prove the result, see e.g.~\cite{mpt,peterzil-starchenko:alw,pila-tsimerman:ax-schanuel}.  

Another well-understood class is when $f$ satisfies a certain kind of differential equation over $\C[z]$, in which case one may appeal to strong results in differential algebra that classify the geometry of the space of solutions of the equation, see \cite{ax,cfn:alw,bcfn:ax-schanuel}.  
However, our proof must take a different approach. 

To begin with, $\Gamma$ does not fall into the desired quotient map setup, and even the restriction of $\Gamma$ to a ``fundamental domain'' is not part of a known o-minimal expansion of $\R$ (see \cite{O-minimalComplexGamma}). 
On the other hand, a famous theorem of H\"older \cite{holder:gamma} shows that $\Gamma$ is a differentially transcendental function, so one cannot appeal to methods from differential algebra.
Our results therefore need new methods. 

Since $\Gamma$ does satisfy some difference equations (like $\Gamma(z+1)=z\Gamma(z)$), we will seek help from difference algebra, mainly in the form of \cite[Th\'eor\`eme 4.2]{Har08}.
However, this will not be enough. 
Most of our methods exploit the ``very transcendental'' nature of $\Gamma$, which for us will be reflected in two ways:
\begin{enumerate}[(a)]
    \item The real analytic curves in $\C$ on which the modulus of $\Gamma(z)$ is constant are not algebraic.
    The same holds for the curves on which the argument of $\Gamma(z)$ is constant. 
    \item Given $c\in\C$, the fiber $\Gamma^{-1}(c)$ is exponentially close to negative integers, which means that for every sufficiently large positive integer $n$, there exists a unique $z_n\in\Gamma^{-1}(c)$ such that the distance between $-n$ and $z_n$ decreases exponentially fast to zero. 
\end{enumerate}

\subsection{Structure of the paper}
After recalling various properties of $\Gamma$ in \S\ref{sec:prelims} and reviewing a few results from difference algebra \S\ref{sec:diffalg}, we study the growth behavior of $\Gamma$ in \S\ref{sec:growth}, especially focusing on its behavior near the poles at negative integers.
In \S\ref{sec:almost-integer-valued} we show that algebraic functions which are ``exponentially close'' to being integer valued, are in fact integer valued. 
The main theorems are proved in \S\ref{sec:mainproof}.

\section{Preliminaries}
\label{sec:prelims}
In this section, we recall some standard facts about the $\Gamma$ function and refer the reader to \cite{remmert} for details. 
The $\Gamma$ function can be defined as a meromorphic function on $\C$, having a simple pole at every non-positive integer, via the following infinite product
\begin{equation}
\label{eq:gammaweierstrass}
    \Gamma(z)\coloneqq  \frac{\exp(-\gamma z)}{z}\prod_{k=1}^{+\infty}\left(1+\frac{z}{k}\right)^{-1}\exp\left(\frac{z}{k}\right),
\end{equation}
where $\gamma\in\R^+$ denotes the Euler--Mascheroni constant 
\begin{equation*}
    \gamma \coloneqq \lim_{n\rightarrow+\infty}\left(\sum_{k=1}^{n}\frac{1}{k}-\log(n)\right).
\end{equation*}
Some basic properties of $\Gamma$ are that $\Gamma(n)=(n-1)!$ for every positive integer $n$, $\Gamma(x)>0$ whenever $x>0$,  $\overline{\Gamma(z)} = \Gamma(\overline{z})$, and $|\Gamma(z)|\leq\Gamma(\Re(z))$ whenever $\Re(z)>0$.

\subsection{Functional equations}
\label{subsec:diffeqs}
 $\Gamma$ is known to satisfy the following functional/difference equations:
\begin{enumerate}[(a)]
    \item For every $z\in\C\setminus\Z_{\leq 0}$,
    \begin{equation}
    \label{eq:diffeq1}
        \Gamma(z+1) = z\Gamma(z).
    \end{equation}
    \item For every $z\in\C$,
    \begin{equation}
    \label{eq:diffeq2}
        \Gamma(z)\Gamma(1-z) = \frac{\pi}{\sin(\pi z)}.
    \end{equation}
    \item For every $z\in\C\setminus\R_{\leq0}$ and for every integer $n\geq 2$,
    \begin{equation}
    \label{eq:diffeq3}
        \prod_{k=0}^{n-1}\Gamma\left(z + \frac{k}{n}\right) = (2\pi)^{\frac{1}{2}(n-1)}n^{\frac{1}{2}-nz}\Gamma(nz).
    \end{equation}
\end{enumerate}  

\begin{defn}
  Given $k\in\N$ we define the \emph{$k$-th Pochhammer polynomial} $m_k(X)\in\Z[X]$ as $m_0(X) = 1$, and for $k>0$
\[m_k(X)\coloneqq X(X+1)\cdots(X+k-1).\]
When $k<0$ we define the \emph{Pochhammer rational function}
$m_k(X)\coloneqq\frac{1}{(X-1)\cdots(X-|k|)}$.
Thus, $\Gamma(z+k) = m_k(z)\Gamma(z)$ for all $k\in\Z$.  
\end{defn}

\subsection{Examples of functional properties of \texorpdfstring{$\Gamma$}{Gamma} not covered by our main results}

In \S\ref{subsec:functrans} we gave an example of a curve in $\C^3$ being mapped under $\Gamma$ to a surface in $\C^3$.
That example was built using \eqref{eq:diffeq1}.
Below we show that \eqref{eq:diffeq2} and \eqref{eq:diffeq3} produce interesting examples of functional properties of $\Gamma$ which are not covered by our main theorems. 
New methods may be needed to extend our results to account for these examples. 
 
\begin{ex} 
\begin{enumerate}[(a)]
    \item We can combine (\ref{eq:diffeq2}) with various trigonometric identities. 
    For example, from $\cos(\pi z)^2+\sin(\pi z)^2 =1$ we obtain
    \[\frac{1}{\Gamma(z)^2\Gamma(1 - z)^2} + \frac{1}{\Gamma\left(\frac{1}{2} + z\right)^2\Gamma\left(\frac{1}{2} - z\right)^2} = \frac{1}{\pi^2}.\]
    Or to put it polynomially:
    \[\Gamma(z)^2\Gamma(1-z)^2 + \Gamma\left(\frac{1}{2} + z\right)^2\Gamma\left(\frac{1}{2} - z\right)^2 = \frac{1}{\pi^2}\Gamma(z)^2\Gamma(1-z)^2\Gamma\left(\frac{1}{2} + z\right)^2\Gamma\left(\frac{1}{2} - z\right)^2.\]
    Here $\Gamma$ takes a curve in $\C^4$ to a hypersurface in $\C^4$.
    
    \item Finally, we use (\ref{eq:diffeq3}) with $n=2$ to obtain
    \[\Gamma(z) \Gamma\left(z + \frac{1}{2}\right) = 2^{1 - 2z}\sqrt{\pi}\Gamma(2z)\quad \mbox{ and }\quad\Gamma(2z) \Gamma\left(2z + \frac{1}{2}\right) = 2^{1 - 4z}\sqrt{\pi}\Gamma(4z).\]
    Squaring the first identity and then dividing by the second identity gives 
    \[\Gamma(z)^2\Gamma\left(z + \frac{1}{2}\right)^2\Gamma(4z) = 2\sqrt{\pi}  \Gamma\left(2z + \frac{1}{2}\right) \Gamma(2z).\]
    Here $\Gamma$ takes a curve in $\C^5$ to a hypersurface in $\C^5$.
\end{enumerate}    
\end{ex}

\subsection{Argument of complex numbers}
When talking about the argument of a complex number $z\in\C$, we will always assume that the argument function $\arg(z)$ is defined on all of $\C^\times$. 
Of course, this implies that $\arg(z)$ is not continuous on all of $\C^\times$, but it is continuous when restricted to the domain of a specific branch of logarithm. 
We express the argument function as $\arg:\C^\times\to(-\pi,\pi]$ and it is given by $\arg(z) = \arccos\left(\frac{\Re(z)}{|z|}\right)$ when $\Im(z)>0$ and by $\arg(z) = -\arccos\left(\frac{\Re(z)}{|z|}\right)$ when $\Im(z)<0$.
We extend the argument function to zero by defining $\arg 0 := 0$.

\subsection{Stirling's formula}
\label{subsec:stirling}
We recall Stirling's formula: there is a holomorphic function $\mu(z)$ on $\C\setminus\R_{<0}$ such that  
\begin{equation}
\label{eq:stirlingeq}
    \Gamma(z) = \sqrt{2\pi}z^{z-\frac{1}{2}}\exp(-z)\exp(\mu(z)),
\end{equation}
and for any $\theta\in(0,\pi)$, if $z$ approaches $\infty$ while satisfying $-\pi+\theta<\arg(z)<\pi-\theta$, then $\mu(z)$ approaches 0.
There are more precise descriptions of $\mu(z)$, see \cite[Chapter 2, \S 4.2]{remmert}. 
In particular, given $\theta \in \left(\frac{\pi}{2},\pi\right)$, there is $M_{\theta}>0$ such that whenever $z$ satisfies $|\arg (z)|<\theta$ and $|z|>M_{\theta}$, we have $|\mu(z)| < 1$, and so $\frac{1}{e}\leq|\exp(\mu(z))|\leq e$, where $e=\exp(1)$. 
Thus, if we write $z=x+iy$ with $x,y\in\R$, we get that when $|\arg (z)|<\theta$ and $|z|> M_{\theta}$
\begin{equation}
\label{eq:stirlingineq}
   \frac{\sqrt{2\pi}}{e}|z|^{x-\frac{1}{2}}\exp\left(- y\arg(z) - x\right) \leq|\Gamma(z)| \leq e\sqrt{2\pi}|z|^{x-\frac{1}{2}}\exp\left(- y\arg(z) - x\right).
\end{equation}

\section{Difference Algebra}
\label{sec:diffalg}
Let $\mathcal{M}(\C)$ be the field  of meromorphic functions on $\C$, and let $K$ be the subfield of $\mathcal{M}(\C)$ of functions that are $1$-periodic, i.e.
\[K = \{f\in\mathcal{M}(\C) : \forall z\in\C (f(z+1) = f(z))\}.\]
As mentioned in the introduction, a famous theorem of H\"older \cite{holder:gamma} says that $\Gamma$ is differentially transcendental over $\C(z)$, and now there are many known proofs of this result, even showing that $\Gamma$ is differentially transcendental over $K(z)$. 
However, $\Gamma$ satisfies the difference equation \eqref{eq:diffeq1}, and some functional properties can be obtained from this through difference algebraic methods. 

We will use the following theorem of Hardouin. 
\begin{thm}[{\cite[Th\'eor\`eme 4.2]{Har08}}]
\label{thm:hardouin}
    Let $a_1, \dots, a_n \in \C(z)^\times$ be non-zero rational functions, and let $f_1, \dots, f_n$ be non-zero meromorphic functions on $\C$ satisfying the equations $f_i(z + 1) = a_i(z)f_i(z)$ for $i=1,\ldots,n$. 
    Then the functions $f_1(z), \ldots, f_n(z)$ are algebraically dependent over $K(z)$ if and only if there exist integers $r_1, \dots, r_n \in \Z$ not all $0$ and a rational function $h \in \C(z)^\times$ such that
    \[a_1(z)^{r_1} \cdots a_n(z)^{r_n} = \frac{h(z + 1)}{h(z)}.\]
\end{thm}

The following corollary is a special case of \cite[Corollaire 4.4]{Har08}, as well as a special case of \cite[Theorem 1]{divizio-pellarin}.

\begin{cor}
\label{cor:holder1}
    For any positive integer $n$, any $a\in\C^\times$ and any $b_1,\ldots,b_n\in\C$, the functions $\Gamma(az+b_1),\ldots,\Gamma(az+b_n)$ are algebraically dependent over $K(z)$ if and only if there are $i,j\in\{1,\ldots,n\}$ distinct such that $b_i-b_j\in\Z$.\footnote{In fact the cited works prove a stronger version of this corollary by replacing \emph{algebraic dependency over $K(z)$} with \emph{differential algebraic dependency over $K(z)$}.}
\end{cor}

\begin{cor}
    \label{cor:twofunctions}
    Given $a\in\C^\times$ and $b\in\C$, the functions $\Gamma(z)$ and $\Gamma(az+b)$ are algebraically dependent over $K(z)$ if and only if $b\in\Z$ and $a=\pm1$. 
    
    Furthermore, $\Gamma(z)$ and $\Gamma(az+b)$ are algebraically dependent over $\C$ if and only if $b=0$ and $a=1$.
\end{cor}
\begin{proof}
    First observe that if $b\in\Z$ and $a=\pm1$, then $\Gamma(z)$ and $\Gamma(az+b)$ are algebraically dependent over $K(z)$ by \eqref{eq:diffeq1} and \eqref{eq:diffeq2} (note that $\sin(\pi(z-b))$ is algebraic over $K$).
    
    Conversely, assume $\Gamma(z)$ and $\Gamma(az+b)$ are algebraically dependent over $K(z)$.
    Shifting by $z\mapsto z+1$ and using the exchange principle of algebraic dependence, we get that $\Gamma(az+b)$ and $\Gamma(az+b+a)$ are algebraically dependent over $K(z)$, giving that $a\in\Z\setminus\{0\}$ by Corollary \ref{cor:holder1}.

    Since $a$ is an integer, we can now apply Theorem \ref{thm:hardouin} to $\Gamma(z)$ and $\Gamma(az+b)$ to get that there are $r_1, r_2 \in \Z$ not both zero and coprime polynomials $ p(z),q(z) \in \C[z]$ such that $z^{r_1}m_a(az+b)^{r_2} = \frac{p(z+1)q(z)}{p(z)q(z+1)}$.
    The rational function on the right hand side has degree 0, hence $r_2\deg(m_a(X)) + r_1=0$, so in particular $r_2\neq 0$. 
    Also, 
    \[\lim_{z\to\infty}\frac{p(z+1)q(z)}{p(z)q(z+1)}=1,\]
    which implies $a^{r_2}=1$, thus $|a|=1$. 
    If $a=1$, Corollary \ref{cor:holder1} shows $b\in\Z$.
    If $a=-1$, then using \eqref{eq:diffeq2} we obtain that $\Gamma(z)$ and $\Gamma(-az-b) = \Gamma(z-b)$ are algebraically dependent over $K(z)$, so again we get $b\in\Z$. 
    This completes the proof of the first part of the corollary.

    For the ``furthermore'' statement, one direction is obvious, so assume that $\Gamma(z)$ and $\Gamma(az+b)$ are algebraically dependent over $\C$, $a=\pm 1$ and $b\in\Z$.
    If $a=-1$, then by \eqref{eq:diffeq2} the functions $\Gamma(z)$ and $\Gamma(z+1-b)\sin(\pi(b-z))$ are algebraically dependent over $\C$. 
    But then by \eqref{eq:diffeq1} we would get that $\Gamma(z)$ and $m_{1-b}(z)\sin(\pi(b-z))$ are algebraically dependent over $\C$, and since $m_{1-b}(z)\sin(\pi(b-z))$ satisfies a differential equation over $K(z)$, then $\Gamma(z)$ would satisfy a differential equation over $K(z)$, contradicting H\"older's theorem.
    So we must have $a=1$, but if $\Gamma(z)$ and $\Gamma(z+b) = m_b(z)\Gamma(z)$ are algebraically dependent over $\C$, then this would imply that $\Gamma(z)$ is algebraic over $\C(z)$ (again a contradiction), unless $b=0$. 
\end{proof}

\section{Growth Behavior of \texorpdfstring{$\Gamma$}{Gamma}}
\label{sec:growth}
We first look at the growth behavior of $\Gamma$ in the upper right quadrant; since $\overline{\Gamma(z)} = \Gamma(\overline{z})$, this suffices to describe the behavior of $\Gamma$ in the right half plane.
By \cite[Proposition 2.5]{EP25}, for each $r \in \R_{> 0}$, there is a $C^1$ function $y_r \colon \R_{>\max\{\alpha,\Gamma^{-1}(r)\}} \to \R_{>0}$ which satisfies $|\Gamma(x + iy_r(x))| = r$ (here $\alpha\approx1.4616\ldots$ denotes the only zero of $\Gamma'$ on $\R_{>0}$, see \cite{uchiyama}). 
It was proved in \cite[Lemma 2.6]{EP25} that $y_r'(x) \ge 2 \log \lfloor x\rfloor - 2$, and hence there exists a constant $c > 0$ such that $y_r(x) \ge 2x \log \left(\frac{x}{c}\right)$ for all large enough $x$. 
We now give a more precise growth rate for $y_r$.

\begin{lem}\label{lem:Growth rate of fixed modulus curves}
    There exist constants $c, c' > 0$ independent of $r$ and $N>0$ depending on $r$ such that $c \log x \le y_r'(x) \le c' \log x$ and $cx\log x \le y_r(x) \le c'x\log x$ for all $x \ge N$.
\end{lem}
\begin{proof}
    The existence of $c$ follows from \cite[Lemma 2.6]{EP25}, and we alter the proof of that lemma to find $c'$.
    Recall that for a differentiable product of differentiable functions $f(x) = \prod_{k=0}^{+\infty}f_{k}(x)$ we have $f'(x) = f(x)\sum_{n=0}^{+\infty}\frac{f_{n}'(x)}{f_{n}(x)}$.
    Applying this to the modulus of (\ref{eq:gammaweierstrass}) and using $|\Gamma(x+iy_r(x))| = r$ gives
    \begin{align*}
        0 &= \frac{d}{dx}\left(\frac{\exp(-\gamma x)}{\sqrt{x^2+y_r(x)^2}} \prod_{n=1}^{+\infty} \frac{\exp(x/n)}{\sqrt{\left(1+\frac{x}{n}\right)^2 + \left(\frac{y_r(x)}{n}\right)^2}}\right) \\
            &= r\left(-\gamma -\frac{y_r(x)y_r'(x)+x}{x^2+y_r(x)^2} + \sum_{n=1}^{+\infty} \left(\frac{1}{n}-\frac{y_r(x)y_r'(x)+n+x}{(n+x)^2+y_r(x)^2}\right)\right).
    \end{align*}
    Solving for $y_r'$ gives
    \begin{equation}\label{eq:modulus implicit derivative}
        y_r'(x) = \frac{-\gamma - \frac{x}{x^2 + y_r(x)^2} + \sum_{n = 1}^{+\infty} \left(\frac{1}{n} - \frac{n + x}{(n + x)^2 + y_r(x)^2}\right)}{\frac{y_r(x)}{x^2 + y_r(x)^2} + \sum_{n = 1}^{+\infty} \frac{y_r(x)}{(n + x)^2 + y_r(x)^2}}.
    \end{equation}
    Addressing the denominator, we have
    \begin{equation}\label{eq:denominator of fixed mod derivative}
        \frac{y_r(x)}{x^2 + y_r(x)^2} + \sum_{n = 1}^{+\infty} \frac{y_r(x)}{(n + x)^2 + y_r(x)^2} \ge \int_{x}^{+\infty} \frac{y_r(x)}{t^2 + y_r(x)^2}dt = \frac{\pi}{2} - \arctan\left(\frac{x}{y_r(x)}\right).
    \end{equation}
    We note that there is a constant $N = N(r)$ such that $y_r(x) \ge x$ for all $x \ge N$.
    This is a consequence of the mean value theorem, using that $y_r'(x) \ge 2 \log \lfloor x \rfloor - 2$ holds for all $x > \max\{\alpha,\Gamma^{-1}(r)\}$, and that $y_r(x)>0$ for all $r>0$ and $x \in \dom(y_r)$.
    So we may assume $x$ is large enough so that $y_r(x) \ge x$ and hence the denominator is bounded below by $\frac{\pi}{4}$.

    For large enough $x$,  (\ref{eq:stirlingineq}) gives 
    \[r = |\Gamma(x + iy_r(x))| \le e\sqrt{2\pi} \exp\left(\left(x - \frac{1}{2}\right)\log \sqrt{x^2 + y_r(x)^2} - y_r(x)\arg(x + iy_r(x)) - x\right).\]
    If $x\geq N$ then $\arg(x + iy_r(x)) \ge \frac{\pi}{4}$, which gives
    \[r \le e\sqrt{2\pi} \exp \left(x \log \left(\sqrt{2}y_r(x)\right) - \frac{\pi}{4}y_r(x)\right).\]
    This implies $x > \sqrt{y_r(x)}$ for all large enough $x$ because if not, i.e., if $x \le \sqrt{y_r(x)}$ for arbitrarily large values of $x$, then we would have
    \begin{align*}
        \liminf_{x\to+\infty}\left(x \log \left(\sqrt{2}y_r(x)\right) - \frac{\pi}{4}y_r(x)\right) \le 
            \liminf_{x\to+\infty}\left(\sqrt{y_r(x)} \log \left(\sqrt{2}y_r(x)\right) - \frac{\pi}{4}y_r(x)\right) = -\infty.
    \end{align*}
    Thus, we may assume $x$ is large enough so that $y_r(x) < x^2$.
    This yields an upper bound on the numerator of $y'_r(x)$:
    \begin{align*}
        \sum_{n = 1}^{+\infty} \left( \frac{1}{n} - \frac{n + x}{(n + x)^2 + y_r(x)^2}\right) \le& \sum_{n = 1}^{+\infty} \left(\frac{1}{n} - \frac{1}{n + x + \frac{x^4}{n + x}}\right)\\
        \le & \sum_{n = 1}^{+\infty} \left( \frac{1}{n} - \frac{1}{n + 2x^3}\right) \\
        \le & \sum_{n = 1}^{\lceil 2x^3\rceil} \frac{1}{n} \\
        \le & 1+ \int_1^{2x^3}\frac{1}{t}dt = 1+\log\left(2x^3\right)
    \end{align*}
    Putting the bounds for numerator and denominator of $y_r'(x)$ together, we get
    \[y_r'(x) \le \frac{1+\log(2)+3\log(x)}{\pi/4}.\]
    We now choose $c'>0$ so that for all $x>\alpha$ we get 
    \[\frac{1+\log(2)+3\log(x)}{\pi/4}\leq c'\log(x).\]

    Now we show the desired bounds on $y_r$ hold.
    Since $y_r'(x)-c\log x \ge 0$, we have 
    \begin{align*}
        0 &\le \int_{\max\{\alpha,\Gamma^{-1}(r)\}}^x(y_r'(t)-c\log t)dt \\
            &=y_r(x)-cx\log x+cx +c\log(\max\{\alpha,\Gamma^{-1}(r)\})- y_r(\max\{\alpha,\Gamma^{-1}(r)\}).
    \end{align*}
    Decreasing $c$ and increasing $N$ as needed gives the desired lower bound on $y_r$.
    An analogous argument gives the upper bound on $y_r$, possibly increasing $N$ once again.
\end{proof}

\subsection{The fibers of \texorpdfstring{$\Gamma$}{Gamma}}
\label{subsec:fibers}

As already alluded to in the introduction, our proof will exploit the ``very transcendental'' nature of $\Gamma$, which in practice means a careful study of the asymptotic behavior of the fibers of $\Gamma$ around negative integers. 

\begin{remark}\label{rmk:fibers of Gamma in vertical strip}
    Stirling's formula (\ref{eq:stirlingineq}) implies $|\Gamma(x+iy)|$ decays exponentially to zero as $|y| \to +\infty$ with $x$ fixed.
    As a consequence, any vertical strip in $\C$ contains at most finitely many elements of $\Gamma^{-1}(c) = \{z\in\C : \Gamma(z) = c\}$, for any $c\in\C^\times$. 
    Moreover, if $C$ is an algebraic curve and $N \in \N$, then there are only finitely many $z$ with $|\Re(z)|\le N$ and $(z,\Gamma(z)) \in C$.
\end{remark}

\begin{lem}
     \label{lem:boundedim}
    There is $M>0$ such that for every $c\in\C^\times$ and every $z\in\Gamma^{-1}(c)$ satisfying $\Re(z)<\frac{1}{2}$, we have $|\Im(z)| \leq \max\{M, -\log |c|\}$.
     \end{lem}
     \begin{proof}
     Write $z = x + iy$ with $x, y \in \R$. Then $|\sin( \pi z)|^2 = \sin(\pi x)^2 + \sinh(\pi y)^2$, which gives $|\sin(\pi z)|\geq |\sinh(\pi y)|$. 
     Also $|\sinh(\pi y)| \ge \frac{\exp(\pi |y|) - 1}{2} \ge \frac{\exp(\pi |y|)}{4}$ for $|y| \ge 1$.
     Let $M = M_{\frac{3\pi}{4}}>0$ be a constant which guarantees (\ref{eq:stirlingineq}) on $\left\{z:|\arg z|<\frac{3\pi}{4}\right\}$.
     The identity  $|\Gamma(z)| =\left|\overline{\Gamma(z)}\right| = |\Gamma(\bar{z})|$
     combined with (\ref{eq:stirlingineq}) and the fact that $\arg(z) = -\arg(\overline{z})$ whenever $z\notin\R_{\leq0}$ proves that for all $x < \frac{1}{2}$ and $|y| > M_{\frac{3\pi}{4}}$ we have
     \begin{align}
         |\Gamma(z)| &= \frac{\pi}{|\Gamma(1-z)||\sin(\pi z)|} \nonumber \\ 
         & \leq \frac{e\pi}{\sqrt{2\pi}|1-z|^{\frac{1}{2}-x}\exp(y\arg(1-z)-(1-x))|\sinh(\pi y)|} \nonumber\\
         \label{eqn:Bound on Gamma size}
         &\le \frac{4\pi e|1 - z|^{x - \frac{1}{2}}\exp(1 - x)}{\sqrt{2\pi}\exp(\pi |y|)\exp(y\arg(1-z))}\\
         &\le \frac{2e\sqrt{2\pi e}\left|\frac{e}{1-x}\right|^{\frac{1}{2} - x}}{\exp\left(\frac{\pi}{2} |y|\right)} \nonumber \\
         &\le \frac{C}{\exp\left(\frac{\pi}{2} |y|\right)}\nonumber
     \end{align}
     for some constant $C > 1$.
     This is strictly smaller than $|c|$ whenever $|y| > \frac{2}{\pi} \log\left( \frac{C}{|c|}\right)$.
     Replacing $M$ by $\max\left\{M, \frac{6}{\pi}\log(C)\right\}$ if necessary gives that if $|y|>\max\{M,-\log|c|\}$ and $|c|<1$, then 
     \[3|y|> M -2\log|c|\geq \frac{6}{\pi} \log(C) - 2\log|c|\geq \frac{6}{\pi}\log\left(\frac{C}{|c|}\right).\]
     Thus we have shown that $\Gamma(z) = c$ has no solutions satisfying $\Re (z) < \frac{1}{2}$ and $|\Im(z)| > \max\{M, -\log |c|\}$.
     \end{proof}

\begin{defn}
    Given a positive real number $x$, let $N_0(x)$ be the smallest positive integer such that
        \[\left|\Gamma\left(\frac{1}{2}-N_0(x)\right)\right| = \left|m_{-N_0(x)}\left(\frac{1}{2}\right)\Gamma\left(\frac{1}{2}\right)\right|=\frac{\sqrt{\pi}}{\left|\frac{1}{2}-1\right|\cdots\left|\frac{1}{2}-N_0(x)\right|} < x,\]
        where $m_{-n}(z)$ denotes a Pochhammer function
\end{defn}

\begin{lem}\label{lem:Left Half Plane Gamma solutions}
    Given $c \in \C^{\times}$ let $N_c$ be the smallest positive integer satisfying
    \[N_c \geq\max\left\{\exp\left(\frac{3}{2}\right), \left\lceil\max\{M,N_0(|c|), -\log |c|, \log |c|\}\right\rceil\right\}\]
    and $\Gamma(N_c) \ge \frac{4\pi}{|c|}$.
    Then for every integer $n \ge N_c$, there exists a unique $z_n \in \Gamma^{-1}(c)$ satisfying $-n-\frac{1}{2} \leq \Re(z_n) \leq -n+\frac{1}{2}$.
    Furthermore, $|z_n + n| \le \frac{2\pi }{|c|\Gamma(n)}$.
\end{lem}
\begin{proof}
     Choose $c\in\C^\times$ and let $M>0$ be given by Lemma \ref{lem:boundedim}.
     Set
     \[N = \lceil\max\{M,N_0(|c|), -\log |c|, \log |c|\}\rceil.\]
     Since $N \ge \max\{M, -\log |c|\}$, there are finitely many solutions to $\Gamma(z) = c$ with $-N + \frac{1}{2} \le \Re(z) \le \frac{1}{2}$ by Lemma \ref{lem:boundedim}.

    \begin{claim}
        \label{claim:uniquefibreinstrip}
        For each $n \ge N$ there exists a unique $z_n\in\Gamma^{-1}(c)$ such that $|\Re(z_n) + n| \le \frac{1}{2}$.
    \end{claim}
     \begin{proof}
    Consider the rectangular path $R_n\subseteq\C$ determined by connecting the four corners $-n\pm \frac{1}{2} \pm iN$ via straight line segments. 
    By Lemma \ref{lem:boundedim}, if $\Gamma^{-1}(c)$ intersects the strip $-n-\frac{1}{2}\leq\Re(z)\leq-n+\frac{1}{2}$, then this intersection occurs in the area bounded by $R_n$. 
    In order to show that there is a unique intersection point in this area, we will show that $|\Gamma(z)|<|c|$ for all $z\in R_n$ and use Rouch\'e's theorem (see e.g.~\cite[Chapter VI, \S1, Theorem 1.6]{lang:complexanalysis}).
     
     Given $z$ in either the top or bottom side of $R_n$, we have $|\Im(z)|=N$, so by (\ref{eqn:Bound on Gamma size}) we already get $|\Gamma(z)| < |c|$. 
    If $z$ is on one of the vertical sides of $R_n$, then $\Re(z) = -n\pm\frac{1}{2}$.
    By \cite[Lemma 29]{O-minimalComplexGamma}, the maps $y\in[0,+\infty) \mapsto \left|\Gamma\left(-n-\frac{1}{2}+iy\right)\right|$ and $y\in[0,+\infty) \mapsto \left|\Gamma\left(-n+\frac{1}{2}+iy\right)\right|$ are injective with nonvanishing derivative and decay exponentially to zero as $y \to +\infty$. 
    The same holds for the maps $y\in[0,+\infty) \mapsto \left|\Gamma\left(-n-\frac{1}{2}-iy\right)\right|$ and $y\in[0,+\infty) \mapsto \left|\Gamma\left(-n+\frac{1}{2}-iy\right)\right|$ since $\Gamma(\overline{z}) = \overline{\Gamma(z)}$.
    This shows that the maximum value of $|\Gamma(z)|$ along the left side of $R_n$ is achieved at $z = -n-\frac{1}{2}$ and the maximum value of $|\Gamma(z)|$ along the right side of $R_n$ is achieved at $z = -n+\frac{1}{2}$. 
    Since $\left|\Gamma\left(-n-\frac{1}{2}\right)\right| = \frac{\left|\Gamma\left(-n+\frac{1}{2}\right)\right|}{\left|-n-\frac{1}{2}\right|}<\left|\Gamma\left(-n+\frac{1}{2}\right)\right|$, the maximum value of $|\Gamma(z)|$ along the vertical sides of $R_n$ is achieved at $z = -n+\frac{1}{2}$.
    On the other hand
    \[\left|\Gamma\left(-n+\frac{1}{2}\right)\right| = \left|m_{-n}\left(\frac{1}{2}\right)\Gamma\left(\frac{1}{2}\right)\right| = \frac{\sqrt{\pi}}{\left|\frac{1}{2}-1\right|\cdots\left|\frac{1}{2}-n\right|}<|c|.\]
    Therefore $\left|\frac{1}{\Gamma(z)}\right| > \frac{1}{|c|}$ for all $z\in R_n$. 
    So by Rouch\'e's Theorem, since $z=-n$ is the only zero of $\frac{1}{\Gamma(z)}$ contained in the region bounded by $R_n$, we conclude that $\frac{1}{\Gamma(z)} - \frac{1}{c}$ has a unique solution $z_n$ in the same region.
    \end{proof}

    It remains to prove the upper bound on $|z_n+n|$. 
    Set $N_c$ to be an integer satisfying $N_c\geq\max\{\exp\left(\frac{3}{2}\right),N\}$ and $\Gamma(N_c) > \frac{4\pi}{|c|}$.
    Let $\epsilon < \frac{1}{2}$ be a positive number and $n\geq N_c$ denote a positive integer.
    Formula (\ref{eq:diffeq2}) gives that for all $\theta\in[0,2\pi]$
    \[\frac{1}{|\Gamma(-n + \epsilon \exp(i\theta))|} = \frac{|\sin (\pi \epsilon \exp(i\theta)) \Gamma(1 + n - \epsilon \exp(i\theta))|}{\pi}.\]
    On the one hand, for every $\theta\in[0,2\pi]$, the slope of the line going through $n$ and $1+n-\epsilon\exp(i\theta)$ lies in $(-1,1)$.
    On the other hand, since $n\geq \exp\left(\frac{3}{2}\right)$, the slope of the constant modulus curve $y_{\Gamma(n)}$ starting at $n$ is always at least $1$ by \cite[Lemma 2.6]{EP25}. 
    So the graph of $y_{\Gamma(n)}$ does not intersect the circle of radius $\epsilon$ centered at $1+n$, which is parametrized by $1+n-\epsilon\exp(i\theta)$ for $\theta \in [0,2\pi]$. 
    The graph of $-y_{\Gamma(n)}$ does not intersect the small circle either, and since $|\Gamma(z)| =\left|\overline{\Gamma(z)}\right| = |\Gamma(\bar{z})|$, the graph of $-y_{\Gamma(n)}$ contains all points $z$ in the lower right quadrant such that $|\Gamma(z)|=n$.
    Thus, $|\Gamma(1 + n - \epsilon \exp(i\theta))| > \Gamma(n)$.
    
    Recall that $\sinh(x) \ge x$ for $x \ge 0$ and $\sin(x) \ge x/2$ for $0 \le x \le \frac{\pi}{2}$.
    Since $\epsilon<\frac{1}{2}$, we get that for all $\theta\in[0,2\pi]$
    \[|\epsilon\pi\cos(\theta)|\leq \frac{\pi}{2}\quad\mbox{ and }\quad |\epsilon\pi\sin(\theta)|\leq \frac{\pi}{2},\]
    and using that $\sin(x)$ and $\sinh(x)$ are odd functions, this gives 
    \[|\sin(\epsilon \pi \exp(i\theta))| = \sqrt{\sin(\epsilon \pi \cos (\theta))^2 + \sinh(\epsilon \pi \sin (\theta))^2} \ge \frac{\epsilon}{2}\]
    by splitting up the circle and using the appropriate bound.
    Combining this all together gives that for all $\theta\in[0,2\pi]$
    \[\frac{1}{|\Gamma(-n + \epsilon \exp(i\theta))|} \ge  \frac{\epsilon \Gamma(n)}{2\pi}.\]
    Choose $\epsilon = \frac{2\pi}{|c|\Gamma(n)}$ (which satisfies $\epsilon<\frac{1}{2}$ by our choice of $N_c$) and let $\mathcal{C}$ denote the circle centered at $-n$ of radius $\epsilon$. 
    Then we have shown that $\left|\frac{1}{\Gamma(z)}\right| > \frac{1}{|c|}$ for all $z\in\mathcal{C}$, so by Rouch\'e's Theorem, since $z=-n$ is the only zero of $\frac{1}{\Gamma(z)}$ contained in the disk bounded by $\mathcal{C}$, we conclude that $\frac{1}{\Gamma(z)} - \frac{1}{c}$ has a unique solution $z_n$ in the same region, hence $|z_n + n| \le \epsilon$ as required.
\end{proof}

A less detailed formulation of the following corollary was briefly mentioned in \cite{EP25}.

\begin{cor}\label{cor:Finitely many points on algebraic curve and Gamma}
    Let $C \subset \C^2$ be an irreducible algebraic curve which is not a horizontal nor a vertical line. 
    Suppose that there is a sequence $\{z_n\}_{n \in \N} \subset \C$ of pairwise distinct elements such that $(z_n, \Gamma(z_n)) \in C$ for all $n$.
    Split each $z_n = x_n + iy_n$ into its real and imaginary components and view each $z_n$ as a point $(x_n, y_n) \in \R^2$.
    Let $V = \overline{\{(x_n, y_n)\}}^{\mathrm{Zar}} \subset \R^2$ denote the $\R$-Zariski closure of these points.
    Then either $V=\R^2$, or there is $N\subseteq\N$ such that
    \begin{enumerate}[(a)]
        \item $\N\setminus N$ is finite,
        \item $y_n=0$ for all $n\in N$,
        \item $x_n<0$ for all $n\in N$, and
        \item there is $a\in\Z$ so that for all $n\in N$ and every $\ell_n\in\N$ satisfying $-\ell_n-\frac{1}{2} \leq x_n\leq -\ell_n+\frac{1}{2}$ we have
        \[|x_{n}+\ell_n| \le \frac{2\pi\left|\ell_n+\frac{1}{2}+\epsilon\right|^a}{\Gamma(\ell_n)}.\]
    \end{enumerate}
\end{cor}
\begin{proof}
    Suppose $\dim V = 1$.
    We first show that $x_n<0$ for all but finitely many $n \in \N$.
    
    If not, then by Remark \ref{rmk:fibers of Gamma in vertical strip} there must be a subsequence of $\{z_n : n \in \N\}$ whose real parts approach $+\infty$.
    Take a further subsequence $\{z_{n_k} : k \in \N\}$ such that $y_{n_k} \ge 0$ for all $k \in \N$ or $y_{n_k} \le 0$ for all $k \in \N$. 
    Since the points $(z_n, \Gamma(z_n))$ lie on an algebraic curve and $\Gamma(z_n) \neq 0$, there exist constants $a, b \in \Z_{>0}$ so that
    \begin{equation}\label{eq:Graph of Gamma intersect algebraic curve}
        \frac{1}{|z_n|^{a}} \le |\Gamma(z_n)| \le |z_n|^{b}
    \end{equation}
    whenever $|z_n| \ge 1$.    
    So for all large enough $k$ we have
    \[\frac{1}{|z_{n_k}|\cdots|z_{n_k}+a-1|} \le \frac{1}{|z_{n_k}|^a}\le |\Gamma(z_{n_k})|\le |z_{n_k}|^b\le |z_{n_k}-1|\cdots|z_{n_k}-b-1|.\]
    The first half of these inequalities gives $1 \le |\Gamma(z_{n_k})||z_{n_k}|\cdots|z_{n_k}+a-1| = |\Gamma(z_{n_k}+a)|$. 
    The second half gives $|\Gamma(z_{n_k}-b-1)|=\frac{|\Gamma(z_{n_k})|}{ |z_{n_k}-1|\cdots|z_{n_k}-b-1|} \le 1$.
    By \cite[Proposition 2.5]{EP25}, the set $X=\{w \in \C : \Re(w)\ge2 \wedge \Im(w)\ge 0 \wedge |\Gamma(w)|=1\}$ is the graph of a function $y_1(x)$ with positive slope and no horizontal or vertical asymptotes. 
    $X$ also contains the point $w=2$ because $\Gamma(2) = (2-1)!=1$. 
    By \cite[Lemma 29]{O-minimalComplexGamma}, the map $y\in[0,+\infty) \mapsto \left|\Gamma\left(x+iy\right)\right|$ is injective with nonvanishing derivative for all $x$, so $X$ divides the quadrant $\{w \in \C : \Re(w)\ge2 \wedge \Im(w)\ge 0\}$ into two simply connected regions: the region below $X$ consists of $w$ such that $|\Gamma(w)|>1$ and the region above $X$ consists of $w$ such that $|\Gamma(w)|<1$.
    Suppose now that $y_{n_k}\ge0$ for all $k \in \N$, so that $z_{n_k}$ lies in the upper right quadrant for all large enough $k$.
    Then for all large enough $k$, $1 \le |\Gamma(z_{n_k}+a)|$ means $z_{n_k}$ lies below $X-a$ (here addition is taking place in $\C$), and $|\Gamma(z_{n_k}-b-1)| \le 1$ means $z_{n_k}$ lies above $X+b+1$.
    
    By Lemma \ref{lem:Growth rate of fixed modulus curves}, there are constants $c,c',N>0$ such that $cx\log x \le y_1(x)\le c'x\log x$ for all $x \ge N$.
    Since $(x_{n_k},y_{n_k}) \in V$, there is an algebraic function $\alpha$ whose domain contains an interval of the form $(r,+\infty)$ such that $\alpha(x_{n_k})=y_{n_k}$ for all $x_{n_k}>r$. 
    Consider $\lim_{x\to+\infty}\frac{\alpha(x)}{x\log x}$.
    The limit exists because $\frac{\alpha(x)}{x\log x}$ is continuous and does not oscillate, and $c(x-b-1)\log(x-b-1)\le \alpha(x)\le c'(x+a)\log(x+a)$ implies that the limit is a real number, which is impossible because no algebraic function is asymptotic to $x\log x$. 
    So there cannot be infinitely many $n$ with $x_n>0$ and $y_n\ge0$.
    If $y_{n_k}\le 0$, then since $\Gamma(\overline{w})=\overline{\Gamma(w)}$, an analogous argument shows there cannot be infinitely many $n$ with $x_n>0$ and $y_n \le 0$.

    So $x_n<0$ for all but finitely many $n$. 
    If there were $L>0$ such that $|y_{n}| \ge L$ for infinitely many $n$, then by \eqref{eqn:Bound on Gamma size}
    \[|\Gamma(z_{n})|\le \frac{2e\sqrt{2\pi e}\left|\frac{e}{1-x_{n}}\right|^{\frac{1}{2}-x_{n}}}{\exp\left(\frac{\pi}{2}L\right)}\]
    which decreases exponentially to zero as $x_{n} \to - \infty$, contradicting \eqref{eq:Graph of Gamma intersect algebraic curve}.
    So for any $L>0$ there are only finitely many $n$ with $y_{n}\ge L$, which gives $\left||z_{n}|-|x_{n}|\right|\le |y_n|\to 0$.
    Hence, given $\epsilon>0$ there is $N_\epsilon\in\N$ so that for all $n>N_\epsilon$ we have $\left||z_{n}|-|x_{n}|\right|< \epsilon$.
    By \eqref{eq:Graph of Gamma intersect algebraic curve}, $b\log|z_{n}|\ge \log|\Gamma(z_{n})|$ and $a\log|z_{n}| \ge-\log|\Gamma(z_{n})|$.
    
    Now let $\ell_n$ be the integer nearest to $|x_{n}|$, if there are two such integers, we let $\ell_n$ be the largest of them.
    Then $\ell_n \geq |x_n|-\frac{1}{2}$. 
    On the other hand, for all $n>N_\epsilon$ we have $|x_n| > |z_n|-\epsilon$, so 
    \begin{equation}
    \label{eq:realdensity}
        \ell_n>|x_n|-\frac{1}{2}>|z_n|-\epsilon-\frac{1}{2}\geq \max\{a,b\}\log|z_{n}| \ge \max\{-\log|\Gamma(z_{n})|,\log|\Gamma(z_{n})|\}
    \end{equation}
     for all large enough $n$.
    Let $N_{|\Gamma(z_{n})|}$ denote the integer defined in Lemma \ref{lem:Left Half Plane Gamma solutions}.
    
    \begin{claim}
    \label{claim:realdensity}
        For all large enough $n$, we have $\ell_n \ge N_{|\Gamma(z_{n})|}$,
    \end{claim}
    \begin{proof}
        By \eqref{eq:realdensity}, it remains to check that for all large enough $n$, $\ell_n\geq N_0(|\Gamma(z_n)|)$ and $\Gamma(\ell_n) \ge \frac{4\pi}{|\Gamma(z_n)|}$. 
        For the first inequality, suppose $n\in\N$ is such that $\ell_n<N_0(|\Gamma(z_n)|)$, then by definition $|\Gamma(z_n)| \leq \left|\Gamma\left(\frac{1}{2}-\ell_n\right)\right|$, which then implies $|z_n|^{-a}\leq \left|\Gamma\left(\frac{1}{2}-\ell_n\right)\right|$. 
        This gives (using \eqref{eq:diffeq2})
        \[\frac{\Gamma\left(\ell_n+\frac{1}{2}\right)}{\pi}=\frac{1}{\left|\Gamma\left(\frac{1}{2}-\ell_n\right)\right|}\leq |z_n|^a,\]
        and so by Stirling's formula \eqref{eq:stirlingeq}
        \[\frac{1}{2}\log(2\pi) + \ell_n\log\left(\ell_n+\frac{1}{2}\right)-\left(\ell_n+\frac{1}{2}\right) + \mu\left(\ell_n+\frac{1}{2}\right) \leq a\log|z_n|,\]
        but by \eqref{eq:realdensity}, this inequality can only be satisfied by finitely many values of $n$. 

        For the second inequality, suppose $n \in \N$ is such that $\Gamma(\ell_n) < \frac{4\pi}{|\Gamma(z_n)|}$. Since $||z_n|-|x_n||\to 0$ as $n \to \infty$ and since $\ell_n$ is the integer nearest to $|x_n|$, $||z_n|-\ell_n|<1$ for all large enough $n$. Combining this with \ref{eq:Graph of Gamma intersect algebraic curve} gives
        \[
            \Gamma(\ell_n) < \frac{4\pi}{|\Gamma(z_n)|} \le 4\pi|z_n|^a < 4\pi(\ell_n+1)^a         
        \]
        which is satisfied by only finitely many $n$, as $\Gamma(x)$ grows exponentially as $x \to \infty$.
    \end{proof}
    
    By Claim \ref{claim:realdensity} and Lemma \ref{lem:Left Half Plane Gamma solutions},
    \[|y_{n}|\le |z_{n}+\ell_n|\le \frac{2\pi}{|\Gamma(z_{n})|\Gamma(\ell_n)}\leq \frac{2\pi|z_n|^a}{\Gamma(\ell_n)} \le \frac{2\pi\left|\ell_n+\frac{1}{2}+\epsilon\right|^a}{\Gamma(\ell_n)}\]
    so if $|y_n|>0$ for infinitely many $n$, then $|y_n|$ would decrease exponentially to zero.
    However, since $(x_n,y_n)$ all lie on the algebraic curve $V$, there are inequalities just as in (\ref{eq:Graph of Gamma intersect algebraic curve}) bounding $|y_n|$ between powers of $|x_n|$. 
    So $|y_n| =0$ for all large enough $n$. \qedhere
\end{proof}

\section{Almost-Integer-Valued Algebraic Functions}
\label{sec:almost-integer-valued}
This section is about proving Proposition \ref{prop:AlmostIntegerValued}. 
In \S\ref{subsec:fibers} we showed that fibers of $\Gamma$ are exponentially close to negative integers. 
In this section we show that an algebraic function sending integers to complex numbers that are exponentially close to integers must be integer-valued. 
When proving the main theorems, the idea is to combine the results of these sections to show the bialgebraic varieties of $\Gamma$ are very rigid.

We first recall a result about algebraic functions, see e.g.~\cite[Corollary 13.16]{eisenbud:commutativealg}. 
Let $P(X,Y)\in\C[X,Y]$ be an irreducible polynomial and suppose that the degree of $P$ in the variable $Y$ is a positive integer $m$. 
    Let $\Omega\subseteq\C$ be an unbounded simply connected domain which does not contain critical points of $P$. 
    Then for any algebraic function $f$ on $\Omega$ satisfying $P(z,f(z))=0$, there is a unique integer $k_0$ and a unique sequence of complex numbers $\{a_k\}_{k=k_0}^{+\infty}$ such that
    \[f(z) = \sum_{k=k_0}^{+\infty}a_kz^{\frac{-k}{m}}\]
    for all $z\in\Omega$.
    This series expansion of $f$ is known as the \emph{Puiseux series of $f$ (on $\Omega$) at $\infty$}, and the rational number $\frac{-k_0}{m}$ is known as the \emph{valuation} or \emph{order} of $f$, and we denote it by $\mathrm{ord}(f)=\frac{-k_0}{m}$. 
    Observe that when $f$ is itself a polynomial in $z$ (so $m=1$ and $a_k=0$ for all $k>0$), then $\mathrm{ord}(f)$ coincides with the degree of $f$ as a polynomial in $z$.

\begin{defn}
    Let $f$ be a holomorphic function defined on a simply connected domain $U\subseteq\C$ containing an interval $(r,+\infty)$ for some $r\in\R$. 
    We say that $f$ is \emph{almost-integer-valued} if there exists $C>0$ such that for all $n\in \N$ there is $\ell_n \in \Z$ with $|f(n) -\ell_n| < C\exp(-n)$.
    
    We say $f$ is \emph{integer-valued} if for all large enough $n\in\N$, $f(n)\in\Z$. 
\end{defn}

Recall that for $z\in\C$ and $n\in\N$ one defines the generalized binomial coefficient as
\[\binom{z}{n}\coloneqq\frac{z(z-1)\cdots(z-n+1)}{n!}.\]
Newton's generalized binomial theorem states that given $a,b \in \C$ with $|a|<1$, we have
    \[(1+a)^b = \sum_{n=0}^{+\infty}\binom{b}{n}a^n.\]

\begin{prop}\label{prop:AlmostIntegerValued}
    Let $f$ be an algebraic function defined on a simply connected domain containing $(r,+\infty)$ for some $r \in \R$. 
    If $f$ is almost-integer-valued, then $f$ is integer-valued and $f\in\Q[X]$.
\end{prop}
\begin{proof}
    For a function $g : (r,+\infty) \to \R$, define the difference operator $\Delta$ by $\Delta(g)\coloneqq g(x+1) - g(x)$.
    Observe that if $g$ and $h$ are almost-integer-valued functions, then so are $\Delta(g)$ and $g+h$. For $n \in \N$, $\Delta^n$ denotes the $n$th iterate of $\Delta$, and a straightforward induction on $n \in \N$ shows
    \begin{equation}\label{eq: iterated difference operator}
        \Delta^ng(x) = \sum_{\ell = 0}^{n}(-1)^{n-\ell}\binom{n}{\ell}g(x+\ell).
    \end{equation}
    
    Let $f$ be an algebraic function defined on $(r,+\infty) \subset \R$ which is almost-integer-valued. 
    We may assume $r \in \N$ by increasing it if necessary.
    Expand $f$ as a Puiseux series in the variable $x$ at $+\infty$ as: 
    \[f(x) = \sum_{k=k_0}^{+\infty} a_kx^{-k/m}\]
    for some $k_0 \in \Z$, some $m \in \N$, and some sequence of complex numbers $\{a_k\}_{k=k_0}^{+\infty}$. 
    We first claim that if $f$ is bounded on $(r,+\infty)$, then it must be constant.
    
    \begin{claim}\label{claim:ConstFunc}
        If $k_0 \ge 0$, then $f\equiv a_{0}$.
    \end{claim}
    \begin{proof}
        If $k_0\ge0$, then $\lim_{x \to +\infty} f(x) = a_{0}$. 
        Choose $C>0$ so that for all $n\in \N$ there is $\ell_n\in\Z$ such that $|f(n)-\ell_n|<C\exp(-n)$. 
        Then the sequence $\ell_n$ must converge to $a_{0}$, which means that $\ell_n=a_{0}$ for all large enough $n$.
        But then $f$ attains the value $a_{0}$ infinitely many times, and since $f$ is algebraic, this implies that $f\equiv a_{0}$.
    \end{proof}

    So, from now on, we assume $k_0 < 0$.
        Observe that $\Delta(f)$ is a Puiseux series of order at most $-\frac{k_0}{m}-1$, hence $\Delta^{-k_0+1}(f)$ is an almost-integer-valued Puiseux series of order strictly less than 0. By Claim \ref{claim:ConstFunc} it must be the constant zero function. 
        Choose a polynomial $q(z)\in\C[z]$ of degree at most $-k_0$ so that $f(n)=q(n)$ for all $n\in\{r,\ldots,-k_0+r\}$. 
        On the one hand, for all $n\geq -k_0+1$, (\ref{eq: iterated difference operator}) gives
    \[\Delta^{n}(f-q)(r)= \sum_{\ell=0}^n(-1)^{n-\ell}\binom{n}{\ell}(f-q)(r+\ell) = \sum_{\ell=-k_0+1}^n(-1)^{n-\ell}\binom{n}{\ell}(f-q)(r+\ell).\]
    On the other hand, whenever $n\geq\deg(q)$ we have $\Delta^n(q)=0$, and we just showed that $\Delta^{-k_0+1}(f)=0$. 
    Therefore for all $n\geq -k_0+1$
    \[0=\sum_{\ell=-k_0+1}^n(-1)^{n-\ell}\binom{n}{\ell}(f-q)(\ell),\]
    proving that $f(n)=q(n)$ for all $n\in\N$, and hence $f=q$. 

    Re-indexing, write $f(x) = a_dx^d+ \cdots + a_0$. 
    We will prove that $a_d,\dots,a_0 \in \Q$ by iterating the following argument.
    Let $g(x) = b_nx^n+\cdots+b_0$ be any almost-integer-valued polynomial, all of whose coefficients are nonzero.
    It can be shown by induction on the degree $n$ that
    \[
        \Delta^{n}\left(g\right) = n!b_{n}.
    \]
    Since $\Delta^{n}(g)$ is almost-integer-valued and constant, $n!b_n$ must be an integer, so $b_n \in \Q$.
    Define $\tilde{g}(x) \coloneqq n!g(x)-n!b_nx^n = n!b_{n-1}x^{n-1}+\cdots + n!b_1x+n!b_0$. 
    Then $\tilde{g}$ is an almost-integer-valued polynomial of degree $n-1$, all its coefficients are nonzero, and its leading coefficient is an integer multiple of $b_{n-1}$.  
    Returning to $f$, we may assume that $a_d,\dots,a_0 \ne 0$ since adding on integer multiples of powers of $x$ preserves being almost-integer-valued. 
    Define $g_0 = f$ and $g_{j+1} = \tilde{g_j}$ for $j<d$.
    Then the leading coefficient of $g_j$ is rational, and it is also an integer multiple of the coefficient $a_{d-j}$ in $f$.
    So $a_{d-j} \in \Q$ for all $j = d,\dots,0$.
    
    Finally, since $f \in \Q[X]$, there is some $\ell \in \Z$ such that $\ell f(n) \in \Z$ for all large enough $n$. 
    So $f(n)$ is rational with denominator dividing $\ell$, and if $f(n) \not \in\Z$ then $|f(n)-m_n| \ge \frac{1}{\ell}$. 
    Since $|f(n)-m_n| < Ce^{-n}$ for all large enough $n$, we must have $f(n) \in \Z$ for all but finitely many $n$.
\end{proof}

\section{Proving the Main Results}
\label{sec:mainproof}

\subsection{The case of curves}
\label{subsec:planecurves}
The main step in proving Theorem \ref{thm:equidimALW} is to first prove the result for curves.
We start by considering plane curves. 

\begin{thm}
    \label{thm:Ax Lindemann for Curves}
    Let $V,W\subset\C^2$ be irreducible curves.
    Suppose $(\Gamma(x),\Gamma(y))\in W$ for all $(x,y)\in V\cap U_\Gamma^2$. 
    Then $V$ is trivially bialgebraic. 
\end{thm}
\begin{proof}
    If some coordinate is constant on $V$, then the result is immediate, so we assume that $V$ is not a horizontal nor a vertical line. 
    Let $W_0\coloneqq\{(z_1,z_2)\in W: z_1z_2\neq0\}$. 
    Since we assume $V$ projects dominantly on each coordinate, $\dim W_0=1$.
    Let $W' \subset \C^2$ be given by the Zariski closure in $\C^2$ of the image of $W_0$ under the inversion map $(z_1, z_2) \mapsto \left(\frac{1}{z_1},\frac{1}{z_2}\right)$.
    Since no coordinate is constant on $W$, the same holds for $W'$, so there are only finitely many numbers $c_1,\ldots,c_k\in\C$ such that $(0,c_1),\dots,(0,c_k)\in W'$.

    Let $\alpha(z)$ be an algebraic function defined on a simply connected domain containing an interval of the form $(-\infty,r)$, for some $r\in\R$, such that $(z,\alpha(z))\in V$ for all $z$ in the domain of $\alpha$. 
    Since $\frac{1}{\Gamma(z)}$ has a zero at every nonpositive integer, for every $n > |r|$ there is $j\in\{1,\ldots,k\}$ such that $\frac{1}{\Gamma(\alpha(-n))} = c_j$. 
    Furthermore, reversing the roles of $x$ and $y$ if necessary, we may and will assume that $\mathrm{ord}(\alpha)\leq 1$. 

    \begin{claim}
        \label{claim:curveatintgers}
        The algebraic function $\alpha(z)$ is actually a polynomial with coefficients in $\Q$.
    \end{claim}
    \begin{proof}
    Suppose first that $k=1$ and $c_1=0$. Then $\alpha(-n)\in-\N$ for all $n >|r|$, or in other words, $\alpha(-x)$ is an algebraic function which is integer-valued. 
    It is well-known (and also it follows from Proposition \ref{prop:AlmostIntegerValued}) that then $\alpha(x)\in\Q[X]$.

    Assume now that some $c_j\neq0$. 
    For each $j\in\{1,\ldots,k\}$ with $c_j\neq 0$, let $N_{c_j}$ be the integer defined in Lemma \ref{lem:Left Half Plane Gamma solutions} and define $\widehat{N} = \max\{N_{c_j} : c_j\neq0\}$. 
    Let $M = \min\{|c_1|,\ldots,|c_k|\}$. Then by Lemma \ref{lem:Left Half Plane Gamma solutions}, for every $n\geq\max\left\{|r|,\widehat{N}\right\}$ there is $\ell_n\in\N$ such that 
    \[|\alpha(-n)+\ell_n| \leq \frac{2\pi}{M\Gamma(n)},\]
    which shows that the algebraic function $\alpha(-x)$ is almost-integer-valued. 
    By Proposition \ref{prop:AlmostIntegerValued} we conclude that $\alpha(x)\in\Q[X]$.  
    \end{proof}
    
    Since $\mathrm{ord}(\alpha)\leq1$, Claim \ref{claim:curveatintgers} shows that there are $a,b\in\Q$ such that $\alpha(z) = az+b$. 
    The hypothesis of the proposition gives that the functions $\Gamma(z)$ and $\Gamma(az+b)$ are algebraically dependent over $\C$, so the proposition now follows from Corollary \ref{cor:twofunctions}. 
\end{proof}

The next corollary follows immediately by taking projections onto pairs of coordinates. 

\begin{cor}
    \label{cor:alwforcurves}
    Suppose that $V \subset \C^n$ is an algebraic curve such that $W \coloneqq \overline{\Gamma(V\cap U_\Gamma^n)}^{\mathrm{Zar}} \subset \C^n$ is also an algebraic curve.
    Then $V$ is trivially bialgebraic.
\end{cor}

\begin{prop}
\label{prop:bialggraph}
    Let  $V\subset\C^2$ be an irreducible algebraic curve and let $W\subset\C^4$ be an irreducible algebraic variety of dimension $2$. 
    Suppose that for all $(z_1,z_2)\in V\cap U_\Gamma^2$ we have $(z_1,z_2,\Gamma(z_1),\Gamma(z_2)) \in W$. 
    Then either $V$ is trivially bialgebraic or there is $b\in\Z$ such that $V$ is defined by the equation $X_1 = X_2+b$. 
\end{prop}
\begin{proof}
    We assume $V$ is not trivially bialgebraic. 
    Then the projection of $W$ onto the last two coordinates is dominant, as otherwise we would be back in  the case of Theorem \ref{thm:Ax Lindemann for Curves}.
    
    Let $W_0 = \{(w_1,w_2,w_3,w_4)\in W\cap\C^4 : w_3w_4\neq0\}$, and let $W'$ be the Zariski closure in $\C^3$ of the image of $W_0$ under the map $(w_1,w_2,w_3,w_4)\mapsto\left(w_2,\frac{1}{w_3},\frac{1}{w_4}\right)$. 
    Since no coordinate is constant on $W$, the same holds for $W'$, hence the set
    \[Z=\{(\omega_2,\omega_3,\omega_4)\in W' : \omega_3=0\}\]
    is a proper Zariski closed subset of dimension at most 1.
    Let $Z_1,\ldots,Z_k$ denote all the irreducible components of $Z$, and for each $j\in\{1,\ldots,k\}$ let $p_j\in\C[X,Y]$ be an irreducible polynomial such that for all $(\omega_2,\omega_4)\in\C^2$ we have $p_j(\omega_2,\omega_4)=0$ if and only if $(\omega_2,0,\omega_4)\in Z_j$. 

    Let $\alpha(z)$ denote an algebraic function defined on a simply connected domain containing an interval of the form $(-\infty, r)$ for some $r\in\R$ and such that $(x,\alpha(x))\in V$ for all $x\in(-\infty, r)$. 
    Reversing the roles of the coordinates in $\C^2$ if necessary, we may assume $\mathrm{ord}(\alpha)\leq 1$. 
    Then for every positive integer $n>r$ we have $(-n,\alpha(-n))\in V$, hence $\left(\alpha(-n),0,\frac{1}{\Gamma(\alpha(-n))}\right)\in W'$, which means that for some $j\in\{1,\ldots,k\}$ we have $p_j\left(\alpha(-n), \frac{1}{\Gamma(\alpha(-n))}\right)=0$.
    By Corollary \ref{cor:Finitely many points on algebraic curve and Gamma}, we conclude that for all positive integers $n>r$, either $\alpha(-n)\in-\N$ or $\alpha(-n)$ is exponentially close to a negative integer, which by Proposition \ref{prop:AlmostIntegerValued} then gives $\alpha(-n)\in-\N$ for all large enough $n$ and $\alpha(z)\in\Q[z]$. 
    Since $\mathrm{ord}(\alpha)\leq 1$, there are $a,b\in\Q$ such that $\alpha(z) = az+b$. 
    Thus the functions $\Gamma(z)$ and $\Gamma(az+b)$ are algebraically dependent over $\C(z)$, so Corollary \ref{cor:twofunctions} finishes the proof. 
\end{proof}

Again by taking projections we get the following corollary. 

\begin{cor}
\label{cor:bialggraph}
    Let  $V\subset\C^n$ be an irreducible algebraic curve.
    Suppose $W\subseteq\C^{2n}$ is an irreducible algebraic variety of dimension $2$ such that for all $\boldsymbol{z}\in V\cap U_\Gamma^n$ we have $(\boldsymbol{z},\Gamma(\boldsymbol{z})) \in W$. 
    Then for every $i,j\in\{1,\ldots,n\}$ at least one of the following holds: 
    \begin{enumerate}[(a)]
        \item the coordinate $X_i$ is constant on $V$, 
        \item the coordinate $X_j$ is constant on $V$, 
        \item $i\neq j$ and there is $b_{ij}\in\Z$ such that $z_i = z_j+b_{ij}$ for all $(z_1,\ldots,z_n)\in V$. 
    \end{enumerate}
\end{cor}

\subsection{The general case}

We now prove the main theorems. 

\begin{proof}[Proof of Theorem \ref{thm:equidimALW}]
    We will assume $\dim V=\dim W$ (the other direction is obvious).
    We proceed by induction on the dimension of $V$.
    The case $\dim V=0$ is trivial, and the case $\dim V = 1$ is covered by Corollary \ref{cor:alwforcurves}, so suppose that $\dim V = d > 1$ and the only bialgebraic varieties of $\Gamma$ of dimension less than $d$ are trivially bialgebraic.
    Let $S \in \{1, 2, \dots, n\}$ be such that the image of the projection of $V$ to the first $S$ coordinates of $\C^n$ is of dimension $d - 1$.
    We denote by $\pi_S \colon \C^n \to \C^S$ the projection onto the first $S$ coordinates.
    Let $V' \coloneqq \overline{\pi_S(V)}^{\mathrm{Zar}}$, and let $U \subset V'$ be a dense Zariski open subset over which the fibers of $\pi_S|_V$ are of dimension $1$.
    
    We claim that $\pi_S(W)$ also has dimension $d - 1$.
    It must be of dimension at least $d - 1$ (because it contains $\Gamma(\pi_S(V))$), so suppose for contradiction that it is of dimension $d$.
    Then, we can find some point $\boldsymbol{c} \coloneqq (c_1, \dots, c_S) \in U$ such that the fiber of $V$ above $\boldsymbol{c}$ is of dimension $1$ but the fiber of $W$ above $\Gamma(\boldsymbol{c})$ is of dimension $0$.
    Since $W \supset \Gamma(V)$, this would give an algebraic curve in $V$ being mapped to a single point in $W$, which cannot be since the $\Gamma$ function is locally injective and thus the Jacobian of the map $\Gamma : U_\Gamma^n \to \C^n$ has full rank.
    
    Thus for every $\boldsymbol{c} \coloneqq (c_1, \dots, c_S) \in \C^S \cap U$, the fiber of $V$ above $\boldsymbol{c}$ is a curve whose image lies in the fiber of $W$ above $\Gamma(\boldsymbol{c})$, which is another algebraic curve.
    Applying Corollary \ref{cor:alwforcurves} says that this fiber of $V$ must be trivially bialgebraic.
    However, since the fibers of $V$ under $\pi_S$ vary continuously on $U$ and have dimension 1, there must exist a non-empty subset $T \subset \{S + 1, \dots, n\}$ such that the fibers of $V$ under $\pi_S$ are all defined by equations of the form $X_i = c_i$ for $i \not \in T$, and $X_j=X_k$ for all $j,k\in T$.
    Let $\pi_{T^c}$ be the projection $\C^n \to \C^{n - |T|}$ to the coordinates not in $T$.
    Then our previous argument says that a Zariski dense subset of fibers of $V$ under the map $\pi_{T^c}$ are the diagonal in $\C^{|T|}$.
    Since $V$ is irreducible, we conclude that $V$ is the product of the diagonal in $\C^{|T|}$ with the projection of $V$ under $\pi_{T^c}$.
    We now use induction hypothesis on $\pi_{T^c}(V)$, which completes the proof.
\end{proof}

\begin{proof}[Proof of Theorem \ref{thm:bialggraph}]
    This proof is similar to the previous one.
    We proceed by induction on the dimension of $V$.
    The case $\dim V=0$ is trivial, and the case $\dim V = 1$ is covered by Corollary \ref{cor:bialggraph}, so suppose that $\dim V = d > 1$.
    Let $S \subset \{1, 2, \dots, n\}$ be a subset such that the projection of $V$ to $\C^S$ has dimension $d-1$.
    We denote by $\pi_S:\C^n\to\C^S$ the projection to those coordinates with index in $S$.
    Let $U \subseteq \overline{\pi_S(V)}^{\mathrm{Zar}}$ be a dense Zariski open subset over which the fibers of $\pi_S|_V$ are of dimension $1$.

    Let $\psi_S:\C^{2n} \to\C^S\times \C^S$ be defined as $(\boldsymbol{x},\boldsymbol{y})\mapsto(\pi_S(\boldsymbol{x}),\pi_S(\boldsymbol{y}))$. 
    For every $\boldsymbol{z}\in \pi_S(V)$ we have $(\boldsymbol{z},\Gamma(\boldsymbol{z}))\in \psi_S(W)$, and since $\dim\pi_S(V) = d-1$ and $\Gamma$ is a transcendental function we conclude $2(d-1)\leq\dim\psi_S(W)\leq 2d$.
    Let $U_0$ be the set of $\boldsymbol{c}\in U$ such that the dimension of the fiber $W_{(\boldsymbol{c},\Gamma(\boldsymbol{c}))}$  of $W$ above $(\boldsymbol{c},\Gamma(\boldsymbol{c}))$ is $\dim W - \dim\psi_S(W)$.
    By the fiber-dimension theorem $U_0$ is Zariski dense in $\overline{\pi_S(V)}^{\mathrm{Zar}}$. 

    First observe that since $\dim V_{\boldsymbol{c}}=1$ and $\boldsymbol{z} \in V_{\boldsymbol{c}}\implies (\boldsymbol{z},\Gamma(\boldsymbol{z}))\in W_{(\boldsymbol{c},\Gamma(\boldsymbol{c}))}$, then $\dim W_{(\boldsymbol{c},\Gamma(\boldsymbol{c}))}\geq 1$.
    Furthermore, since $\Gamma$ is a transcendental function, we in fact get $\dim W_{(\boldsymbol{c},\Gamma(\boldsymbol{c}))}\geq 2$.
    Therefore $\dim\psi_S(W)=2d-2$ and $\dim W_{(\boldsymbol{c},\Gamma(\boldsymbol{c}))} = 2$.
    
    We can apply Corollary \ref{cor:bialggraph} to $V_{\boldsymbol{c}}$ and $W_{(\boldsymbol{c},\Gamma(\boldsymbol{c}))}$.
    Since the fibers of $V$ under $\pi_S$ vary continuously on $U$ and have dimension 1, there exist $b\in\Z$ and $j,k\in \{1,\ldots,n\}\setminus S$ distinct such that for all $\boldsymbol{c}\in U$ the fiber $V_{\boldsymbol{c}}$ is contained in the hypersurface defined by $X_j=X_k + b$. 
    Let $B\subseteq\{1,\ldots,n\}$ be the set defined by the property $j\in B$ if and only if there are $k\in\{1,\ldots,n\}$ and $b\in \Z$ such that $j\neq k$ and $V$ is contained in the hypersurface defined by $X_j=X_k + b$. 
    The above then shows that $V$ is non-empty. 
    
    If $B = \{1,\ldots,n\}$, then we are done.
    Otherwise, let $\pi_{B^c}$ be the projection $\C^n \to \C^{n - |B|}$ to the coordinates not in $B$.
    By construction, a Zariski dense subset of fibers of $V$ under the map $\pi_{B^c}$ are of the desired form in $\C^{B}$.
    Since $V$ is irreducible, we conclude that $V$ is the product of the image in $\C^{B}$ with the projection of $V$ under $\pi_{B^c}$.
    We now use induction hypothesis on $\pi_{B^c}(V)$, which completes the proof.
\end{proof}

\bibliographystyle{alphaurl}
\bibliography{references}

\end{document}